\documentclass[12pt]{article}
\usepackage{amsmath,amsfonts,amscd,amssymb,amsthm,graphicx}
\title{Vector Fields on Product Manifolds}
\author{Stefan Kurz\\
\small Tampere University of Technology\\[-0.8ex]
\small Department of Electronics\\[-0.8ex]
\small Electromagnetics\\[-0.8ex]
\small 33101 Tampere, Finland\\[-0.8ex]
\small \texttt{stefan.kurz@tut.fi}
}

\date{\today~ --~ V1}
\theoremstyle{plain}
\newtheorem*{thm}{Theorem}

\theoremstyle{definition}
\newtheorem*{defn}{Definition}

\theoremstyle{remark}
\newtheorem*{rem}{Remark}

\begin{document}
\maketitle
\begin{abstract}
This short report establishes some basic properties of smooth vector fields on product manifolds. The main results are: (i) On a product manifold there always exists a direct sum decomposition into horizontal and vertical vector fields. (ii) Horizontal and vertical vector fields are naturally isomorphic to smooth families of vector fields defined on the factors.\par
Vector fields are regarded as derivations of the algebra of smooth functions. Basic ideas are taken from Chapter 0 of Ref.\ \cite{Vinogradov2009}.
\end{abstract}
\section*{Basic Properties}
\cite[0.2.25]{Vinogradov2009}
Let $M,N$ be differentiable manifolds. A {\em product manifold} will be denoted by $V=M\times N$. The smooth maps
\[
\pi_M:V\to M\quad\text{and}\quad\pi_N:V\to N
\]
will be called the {\em projection maps} of $V$.\par\bigskip\noindent
\cite[0.2.26]{Vinogradov2009}
For each $n_0\in N$ the smooth map
\[
i_{n_0}:M\to V,\quad m\mapsto(m,n_0)
\]
will be called the {\em embedding at} $n_0$ {\em into} $V$. Similarly, the smooth map
\[
j_{m_0}:N\to V,\quad m\mapsto(m_0,n),
\]
$m_0\in M$ will be called the {\em embedding at} $m_0$ {\em into} $V$.\par\bigskip\noindent
\cite[0.2.9]{Vinogradov2009}
The {\em algebra of smooth functions} on $M$ is denoted by $C^{\infty}(M)$.\par\bigskip\noindent
\cite[0.4.1]{Vinogradov2009}
A {\em smooth vector field} on $M$ is a derivation of $C^\infty(M)$. The $C^\infty(M)$-module of all vector fields on $M$ will be denoted by $\mathfrak{X}(M)$.\par\bigskip\noindent
Let $X(M)$ be some space defined over $M$. Consider the family $x=\{x_n\in X(M)\}_{n\in N}$. The space of families $x$ on $M$ will be denoted by $X(M,N)$, subject to some smoothness requirements.\par
In particular, each $g=\{g_n\in C^\infty(M)\}_{n\in N}$ defines a function
\begin{equation}\label{cinfty1}
f:V\to\mathbb{R},\quad v\mapsto g_n(m),\quad (m,n)=(\pi_M,\pi_N)v.
\end{equation}
If $f\in C^\infty(V)$ the family $g$ is said to be a {\em smooth family of functions}. The space of smooth families of functions is denoted by $C^\infty(M,N)$. Conversely, each $f\in C^\infty(V)$ defines $g\in C^\infty(M,N)$ by
\begin{equation}\label{cinfty2}
g=\{i_n^*f\},
\end{equation}
and the two constructions are inverse to each other. Equations \eqref{cinfty1}, \eqref{cinfty2} define a {\em canonical isomorphism}
\begin{equation}\label{canoniso1}
C^\infty(M,N)\cong C^\infty(V).
\end{equation}
\par\bigskip\noindent
\cite[0.6.10]{Vinogradov2009}
Consider the family $\mathbf{w}=\{\mathbf{w}_n\in\mathfrak{X}(M)\}_{n\in N}$. The family will be called a {\em smooth family of vector fields} on $M$ if it is a derivation on $C^\infty(M,N)$, defined by
\[
\mathbf{w}:\{g_n\}\mapsto\{\mathbf{w}_n(g_n)\}.
\]
The space of smooth families of vector fields will be denoted $\mathfrak{X}(M,N)$.
\begin{rem}
In case factor $N$ is 1-dimensional, we will call $\{\mathbf{w}_n\}$ also a {\em smooth $n$-dependent vector field}.
\end{rem}
\par\bigskip\noindent
Isomorphism \eqref{canoniso1} provides a natural embedding
\begin{equation*}
\iota_{\mathfrak{X}(M,N)}:\mathfrak{X}(M,N)\hookrightarrow\mathfrak{X}(V),\quad\{\mathbf{w}_n\}\mapsto\mathbf{v},
\end{equation*}
where
\begin{equation}\label{injmap}
i_n^*\mathbf{v}(f)=\mathbf{w}_n(i_n^*f),\quad\forall f\in C^\infty(V).
\end{equation}
\par\bigskip\noindent
Moreover, define the $C^\infty(V)$-linear map
\begin{equation*}
\pi_{\mathfrak{X}(M,N)}:\mathfrak{X}(V)\to\mathfrak{X}(M,N),\quad\mathbf{v}\mapsto\{\mathbf{w}_n\},
\end{equation*}
where
\begin{equation}\label{surmap}
\mathbf{w}_n(g)=i_n^*\mathbf{v}(\pi_M^*g),\quad\forall g\in C^{\infty}(M).
\end{equation}
\par\bigskip\noindent
By interchanging the roles of $M$ and $N$ we arrive at similar maps $\iota_{\mathfrak{X}(N,M)}$ and $\pi_{\mathfrak{X}(N,M)}$, respectively.
\begin{defn}\label{defn:horver}
Define the spaces
\begin{subequations}\label{defhorver}
\begin{alignat}{4}
\mathfrak{X}_N(V)&=\{\mathbf{v}\in\mathfrak{X}(V)\,|\,\mathbf{v}(\pi_N^*h)&&=0\quad\forall h&&\in C^{\infty}(N)&&\},\label{defhor}\\[0.3\baselineskip]
\mathfrak{X}_M(V)&=\{\mathbf{v}\in\mathfrak{X}(V)\,|\,\mathbf{v}(\pi_M^*g)&&=0\quad\forall g&&\in C^{\infty}(M)&&\}.\label{defver}
\end{alignat}
\end{subequations}
Vector fields in $\mathfrak{X}_N(V)$ and $\mathfrak{X}_M(V)$ are called {\em horizontal} and {\em vertical vector fields} with respect to the first factor, respectively.
\end{defn}
\begin{rem}
There are equivalent definitions
\begin{alignat*}{4}
\mathfrak{X}_N(V)&=\{\mathbf{v}\in\mathfrak{X}(V)\,|\,(\pi_N^*\boldsymbol{\eta})(\mathbf{v})&&=0\quad\forall\boldsymbol{\eta}&&\in \mathcal{F}^1(N)&&\},\\[0.3\baselineskip]
\mathfrak{X}_M(V)&=\{\mathbf{v}\in\mathfrak{X}(V)\,|\,(\pi_M^*\boldsymbol{\omega})(\mathbf{v})&&=0\quad\forall\boldsymbol{\omega}&&\in \mathcal{F}^1(M)&&\},
\end{alignat*}
where $\mathcal{F}^1$ denotes the space of differential 1-forms.
\end{rem}
\begin{thm}Decomposition of vector fields on product manifolds.\par
\begin{enumerate}
\item \cite[0.4.20]{Vinogradov2009} Every vector field on a product manifold may be decomposed into a horizontal and a vertical component,
\begin{equation}\label{hvdecomp}
\mathfrak{X}(V)=\mathfrak{X}_N(V)\oplus \mathfrak{X}_M(V).
\end{equation}
\item The projection maps related to decomposition \eqref{hvdecomp} are
\begin{equation}\label{projmaps}
\left.\begin{alignedat}{3}
\pi_{\mathfrak{X}_N(V)}&=\iota_{\mathfrak{X}(M,N)}\circ\pi_{\mathfrak{X}(M,N)}:\mathfrak{X}(V)\to\mathfrak{X}_N(V)\\[0.2\baselineskip]
\pi_{\mathfrak{X}_M(V)}&=\iota_{\mathfrak{X}(N,M)}\circ\pi_{\mathfrak{X}(N,M)}:\mathfrak{X}(V)\to\mathfrak{X}_M(V)
\end{alignedat}\quad\right\}.
\end{equation}
\item There are natural isomorphisms
\begin{equation}\label{canoniso2}
\left.
\begin{alignedat}{3}
\iota_{\mathfrak{X}(M,N)}&:\mathfrak{X}(M,N)&&\xrightarrow{\,\,\sim\,\,}\mathfrak{X}_N(V)\\[0.2\baselineskip]
\iota_{\mathfrak{X}(N,M)}&:\mathfrak{X}(N,M)&&\xrightarrow{\,\,\sim\,\,}\mathfrak{X}_M(V)
\end{alignedat}
\quad\right\}.
\end{equation}
\item The following sequence is exact:
\begin{equation}\label{exactseq}
\begin{CD}
0\rightarrow\mathfrak{X}(M,N)@>\displaystyle\iota_{\mathfrak{X}(M,N)}>>\mathfrak{X}(V)@>\displaystyle\pi_{\mathfrak{X}(N,M)}>>\mathfrak{X}(N,M)\rightarrow 0.
\end{CD}
\end{equation}
\end{enumerate}
\end{thm}
\begin{proof}
Let $\{\tilde{\mathbf{w}}_n\}=\pi_{\mathfrak{X}(M,N)}\mathbf{v}$, $\mathbf{v}=\iota_{\mathfrak{X}(M,N)}\{\mathbf{w}_n\}$ in the sequel. Then $\forall g\in C^{\infty}(M)$
\begin{gather}
\tilde{\mathbf{w}}_n(g)\stackrel{\eqref{surmap}}{=}i_n^*\mathbf{v}(\pi_M^*g)
\stackrel{\eqref{injmap}}{=}\mathbf{w}_n(i_n^*\pi_M^*g)=\mathbf{w}_n(g)\nonumber\\[0.3\baselineskip]
\Leftrightarrow\tilde{\mathbf{w}}_n=\mathbf{w}_n\Leftrightarrow\{\tilde{\mathbf{w}}_n\}=\{\mathbf{w}_n\}\nonumber\\[0.3\baselineskip]
\Leftrightarrow\pi_{\mathfrak{X}(M,N)}\circ\iota_{\mathfrak{X}(M,N)}=\mathrm{Id}_{\mathfrak{X}(M,N)},\label{proof1}
\end{gather}
where we took into account $\pi_M\circ i_n=\mathrm{Id}_M$. Equation \eqref{proof1} implies 
\begin{align}
\mathrm{Ker}(\iota_{\mathfrak{X}(M,N)})&=\mathbf{0},\label{keriota}\\
\mathrm{Im}(\pi_{\mathfrak{X}(M,N)})&=\mathfrak{X}(M,N).\label{impi}
\end{align}
Moreover, it follows that $\pi_{\mathfrak{X}_M(V)}$ is idempotent,
\begin{align}
\pi_{\mathfrak{X}_M(V)}\circ\pi_{\mathfrak{X}_M(V)}&=\iota_{\mathfrak{X}(M,N)}\circ\pi_{\mathfrak{X}(M,N)}\circ\iota_{\mathfrak{X}(M,N)}\circ\pi_{\mathfrak{X}(M,N)}\nonumber\\
&=\iota_{\mathfrak{X}(M,N)}\circ\pi_{\mathfrak{X}(M,N)}=\pi_{\mathfrak{X}_M(V)},\label{idempot}
\end{align}
hence a projection. It induces a direct sum decomposition
\begin{alignat}{3}
\mathfrak{X}(V)&=\mathrm{Im}(\pi_{\mathfrak{X}_M(V)})&&\oplus\mathrm{Ker}(\pi_{\mathfrak{X}_M(V)})\nonumber\\
&=\mathrm{Im}(\iota_{\mathfrak{X}(M,N)})&&\oplus\mathrm{Ker}(\pi_{\mathfrak{X}(M,N)}).\label{proof2}
\end{alignat}
Consider $\forall g\in C^\infty(M)$
\begin{gather*}
\{\tilde{\mathbf{w}}_n\}=0\Leftrightarrow \tilde{\mathbf{w}}_n=\mathbf{0}\Leftrightarrow\tilde{\mathbf{w}}_n(g)=0\\[0.3\baselineskip]
\stackrel{\eqref{surmap}}{\Leftrightarrow}i_n^*\mathbf{v}(\pi_M^*g)=0
\stackrel{\eqref{canoniso1}}{\Leftrightarrow}\mathbf{v}(\pi_M^*g)=0,
\end{gather*}
therefore
\begin{align}
\mathrm{Ker}(\pi_{\mathfrak{X}(M,N)})
&=\{\mathbf{v}\in\mathfrak{X}(V)\,|\,\{\tilde{\mathbf{w}}_n\}=\pi_{\mathfrak{X}(M,N)}\mathbf{v}=0\}\nonumber\\[0.3\baselineskip]
&=\{\mathbf{v}\in\mathfrak{X}(V)\,|\,\mathbf{v}(\pi_M^*g)=0\quad\forall g\in C^{\infty}(M)\}\nonumber\\[0.1\baselineskip]
&\stackrel{\makebox[0pt]{\scriptsize\eqref{defver}}}{=}\,\mathfrak{X}_M(V).\label{kerpi}
\end{align}
Pick $f$ in \eqref{injmap} according to $f=\pi_N^*h,h\in C^\infty(N)$, which yields
\[
i_n^*\mathbf{v}(\pi_N^*h)\stackrel{\eqref{injmap}}{=}\mathbf{w}_n(i_n^*\pi_N^*h)=0,
\]
where we took into account $\pi_N\circ i_{n_0}:M\to N$, $m\mapsto n_0$. But then $\forall h\in C^\infty(N)$
\[
i_n^*\mathbf{v}(\pi_N^*h)=0\stackrel{\eqref{canoniso1}}{\Leftrightarrow}\mathbf{v}(\pi_N^*h)=0
\stackrel{\eqref{defhor}}{\Leftrightarrow}\mathbf{v}\in\mathfrak{X}_N(V),
\]
from which we infer that $\mathrm{Im}(\iota_{\mathfrak{X}(M,N)})\subseteq\mathfrak{X}_N(V)\stackrel{\eqref{kerpi}}{=}\mathrm{Ker}(\pi_{\mathfrak{X}(N,M)})$. The direct sum decomposition
\begin{equation}\label{imiota1}
\mathrm{Ker}(\pi_{\mathfrak{X}(N,M)})=\mathrm{Im}(\iota_{\mathfrak{X}(M,N)})\oplus\mathfrak{Z}
\end{equation}
of submodules follows. Equations \eqref{keriota} and \eqref{impi} imply $\dim\mathrm{Im}(\iota_{\mathfrak{X}(M,N)})=\dim M$ and $\dim\mathrm{Ker}(\pi_{\mathfrak{X}(N,M)})=\dim V-\dim N$, respectively\footnote{By \cite[Thm.\ 11.32]{Nestruev2003}, the modules considered here are isomorphic to modules of sections of smooth vector bundles. Therefore, we can identify ranks of modules with dimensions of vector bundles.}. Therefore, we receive from \eqref{imiota1}
\[
\dim V= \dim M+\dim N+\dim\mathfrak{Z},
\]
hence $\mathfrak{Z}=\{\mathbf{0}\}$ and
\begin{equation}\label{imiota}
\mathrm{Im}(\iota_{\mathfrak{X}(M,N)})=\mathfrak{X}_N(V).
\end{equation}
Now \eqref{hvdecomp} follows from \eqref{proof2} with \eqref{kerpi} and \eqref{imiota}, \eqref{projmaps} follows from \eqref{idempot} with \eqref{imiota}, \eqref{canoniso2} follows from \eqref{keriota} and \eqref{imiota}, \eqref{exactseq} follows from \eqref{keriota}, \eqref{kerpi} with \eqref{imiota} and \eqref{impi}.
\end{proof}
\begin{rem}
The Theorem naturally generalizes to product manifolds with more than two factors. Suppose we have $V=M\times N\times L$. Define on top of \eqref{defhorver}
\[
\mathfrak{X}_L(V)=\{\mathbf{v}\in\mathfrak{X}(V)\,|\,\mathbf{v}(\pi_L^*f)=0\quad\forall f\in C^{\infty}(L)\}.
\]
We have the horizontal subspace $\mathfrak{X}_{NL}(V)=\mathfrak{X}_N(V)\cap\mathfrak{X}_L(V)$, and the two vertical subspaces $\mathfrak{X}_{LM}(V)=\mathfrak{X}_L(V)\cap\mathfrak{X}_M(V)$, and $\mathfrak{X}_{MN}(V)=\mathfrak{X}_M(V)\cap\mathfrak{X}_N(V)$, respectively. Decomposition \eqref{hvdecomp} now reads
\[
\mathfrak{X}(V)=\mathfrak{X}_{NL}(V)\oplus\mathfrak{X}_{LM}(V)\oplus\mathfrak{X}_{MN}(V),
\]
with projection map
\[
\pi_{\mathfrak{X}_{NL}(V)}=\iota_{\mathfrak{X}(M,N\times L)}\circ\pi_{\mathfrak{X}(M,N\times L)}:\mathfrak{X}(V)\to\mathfrak{X}_{NL}(V),
\]
and natural isomorphism
\[
\iota_{\mathfrak{X}(M,N\times L)}:\mathfrak{X}(M,N\times L)\xrightarrow{\,\,\sim\,\,}\mathfrak{X}_{NL}(V).
\]
The remaining two projection maps and natural isomorphisms are obtained by cyclic permutation of $M$, $N$, and $L$.\par
There is no exact sequence similar to \eqref{exactseq}, though. It rather holds that $\mathrm{Ker}(\pi_{\mathfrak{X}(M,N\times L)})=\mathfrak{X}_M(V)$, and therefore
\[
\mathrm{Im}(\iota_{\mathfrak{X}(M,N\times L)})=\mathrm{Ker}(\pi_{\mathfrak{X}(N,L\times M)})\cap\mathrm{Ker}(\pi_{\mathfrak{X}(L,M\times N)}).
\]
\end{rem}
\bibliographystyle{plain}
%

%
\end{document}